\documentclass[12pt,twoside]{article}

\usepackage{amsfonts}
\usepackage{amssymb}
\usepackage{amsmath}
\usepackage{graphicx}

\newcommand{\ignore}[1]{}

\def\@begintheorem#1#2{\par\bgroup{\sc #1\ #2. }\it\ignorespaces}
\def\@opargbegintheorem#1#2#3{\par\bgroup{\sc #1\ #2\ (#3). } \it\ignorespaces}
\def\@endtheorem{\egroup}
\newtheorem{theorem}{Theorem}[section]
\newtheorem{corollary}[theorem]{Corollary}
\newtheorem{lemma}[theorem]{Lemma}
\newtheorem{example}[theorem]{Example}
\newtheorem{proposition}[theorem]{Proposition}
\newtheorem{definition}[theorem]{Definition}
\newcommand{\bt}[1]{\begin{theorem}\label{#1}}
\newcommand{\bc}[1]{\begin{corollary}\label{#1}}
\newcommand{\bl}[1]{\begin{lemma}\label{#1}}
\newcommand{\be}[1]{\begin{example}\label{#1}}
\newcommand{\bp}[1]{\begin{proposition}\label{#1}}
\newcommand{\ba}[1]{\begin{algorithm}\rm\label{#1}}
\newcommand{\bd}[1]{\begin{definition}\rm\label{#1}}

\newcommand{\et}{\end{theorem}}
\newcommand{\ec}{\end{corollary}}
\newcommand{\el}{\end{lemma}}
\newcommand{\ee}{\end{example}}
\newcommand{\ep}{\end{proposition}}
\newcommand{\ed}{\end{definition}}
\newcommand{\epr}{{\ \vbox{\hrule\hbox{%
\vrule height1.3ex\hskip0.8ex\vrule}\hrule}}\\\par}

\def\Z{\mathbb{Z}}
\def \G {{\cal G}}

\def \supp {{\rm supp}}

\begin{document}

\title{\bf The Huge Multiway Table Problem}
\author{Shmuel Onn}

\date{\small Technion - Israel Institute of Technology, Haifa, Israel
\\onn@ie.technion.ac.il}

\maketitle

\begin{abstract}
Deciding the existence of an $l\times m\times n$ integer threeway table
with given line-sums is NP-complete already for fixed $l=3$, but is in P
with both $l,m$ fixed. Here we consider {\em huge} tables,
where the variable dimension $n$ is encoded in {\em binary}.
Combining recent results on integer cones and Graver bases,
we show that if the number of {\em layer types} is fixed, then the problem is
in P, whereas if it is variable, then the problem is in NP intersect coNP. Our treatment
goes through the more general class of $n$-fold integer programming problems.
\end{abstract}

\section{Introduction}

Consider the problem of deciding if there is a threeway table with given line-sums:
$$\mbox{is}\ \ \{x\in\Z_+^{l\times  m\times n}\ :\ \sum_i x_{i,j,k}=e_{j,k}
\,,\ \sum_j x_{i,j,k}=f_{i,k}\,,\ \sum_k x_{i,j,k}=g_{i,j}\}\ \ \mbox{nonempty ?}
$$
It is NP-complete already for $l=3$, see \cite{DO1}. But, when both $l,m$ are
fixed, it is decidable in polynomial time \cite{DHOW}, and in fact, in time
which is cubic in $n$ and linear in the binary encoding of $e,f,g$, see \cite{HOR}.
Assume throughout then that $l,m$ are fixed, and regard each table as a
tuple $x=(x^1,\dots,x^n)$ consisting of $n$ many $l\times m$ {\em layers}.
We call the problem {\em huge} if the variable number $n$ of layers is encoded
in {\em binary}. We are then given $t$ {\em types} of layers, where each type
$k$ has its column-sums vector $e^k\in\Z_+^m$ and row-sums vector $f^k\in\Z_+^l$.
In addition, we are given positive integers $n_1,\dots,n_t,n$ with $n_1+\cdots+n_t=n$,
all encoded in binary. A feasible table $x=(x^1,\dots,x^n)$ then must have first
$n_1$ layers of type $1$, next $n_2$ layers of type $2$, and so on, with last
$n_t$ layers of type $t$. The special case of $t=1$ is the case of {\em symmetric}
tables, where all layers have the same row and column sums, and the standard
(non-huge) table problem occurs as the special case of $t=n$ and $n_1=\cdots=n_t=1$.

The huge table problem can be formally defined as follows.

\vskip.2cm\noindent
{\bf Huge Table Problem.} Given $t$ types, line-sums $g\in\Z_+^{l\times m}$,
column-sums $e^k\in\Z_+^m$ and row-sums $f^k\in\Z_+^l$ for $k=1,\dots,t$, and
positive integers $n_1,\dots,n_t,n$ with $n_1+\cdots+n_t=n$, all but $t$
encoded in binary, decide the existence of a feasible table.

\vskip.2cm
Note that while $t$ may be small or fixed, the set of possible layers of type $k$ is
$$
\left\{x^k\in\Z_+^{l\times m}\ :\ \sum_i x^k_{i,j}=e^k_j
\,,\ \sum_j x^k_{i,j}=f^k_i\right\}\ ,
$$
and may have cardinality exponential in the binary encoding of $e^k,f^k$, so it
is not off hand clear how to even write down a single table. But we have the following.
\bt{ThreeWayLineSum}
The huge threeway table problem with all data but $t$ binary encoded is in P
for every fixed $t$ and in NP intersect coNP for $t$ variable and unary encoded.
\et
Since problems in NP intersect coNP are often in P, it is particularly interesting
to know whether the table problem with variable $t$ is decidable in polynomial time.
This is open even for the smallest nontrivial case of $3\times 3\times n$ tables, 
where the input consists of $7t+9$ binary encoded nonnegative integers:
$3t$ row-sums, $3t$ column-sums, $9$ vertical line-sums, and $n_1,\dots,n_t$.
(In fact, removing obvious redundancies among the line-sums, the input amounts to 
$6t+4$ binary encoded nonnegative integers.)

\vskip.5cm
These results follow from broader results which we proceed to describe.
The class of {\em $n$-fold integer programming} problems is defined as follows.
Let $A$ be an $(r,s)\times d$ {\em bimatrix}, by which we mean
a matrix having an $r\times d$ block $A_1$ and $s\times d$ block $A_2$,
$$A\quad=\quad
\left(
\begin{array}{c}
  A_1 \\
  A_2 \\
\end{array}
\right)\quad .
$$
Its {\em $n$-fold product} is the following $(r+sn)\times(dn)$ matrix,
$$A^{(n)}\quad:=\quad
\left(
\begin{array}{cccc}
  A_1    & A_1    & \cdots & A_1    \\
  A_2    & 0      & \cdots & 0      \\
  0      & A_2    & \cdots & 0      \\
  \vdots & \vdots & \ddots & \vdots \\
  0      & 0      & \cdots & A_2    \\
\end{array}
\right)\quad .
$$
The $n$-fold integer programming problem is then the following,
$$\min\left\{wx\ :\ x\in\Z^{dn}\,,\ A^{(n)}x=b\,,\ l\leq x\leq u\right\}\ ,$$
where $w\in\Z^{dn}$, $b\in\Z^{r+sn}$, and $l,u\in\Z_{\infty}^{dn}$
with $\Z_{\infty}:=\Z\uplus\{\pm\infty\}$.
For instance, optimization over multiway tables is an $n$-fold program,
as explained later. By solving the optimization problem we mean,
as usual, either finding an optimal solution, or asserting that the problem
is infeasible, or that the objective is unbounded. In particular,
the optimization problem includes as a special case the feasibility problem.

It was shown in \cite{DHOW,HOW}, building on \cite{AT,HS,SS},
that $n$-fold integer programming for {\em fixed} bimatrix
$A$ can be solved in polynomial time. More recently, in \cite{HOR},
it was shown that for fixed $A$ it can be solved in time which is cubic in $n$ and
linear in the binary encoding of $w,b,l,u$, and that if only the dimensions
$r,s,d$ of $A$ are fixed but $A$ is part of the input, then it can be
solved in time cubic in $n$, polynomial in the unary encoding of $A$,
and linear in the binary encoding of $w,b,l,u$. See \cite{Onn} for a detailed treatment of
the theory and applications of $n$-fold integer programming.

The vector ingredients of an $n$-fold integer program are naturally arranged
in {\em bricks}, where $w=(w^1,\dots,w^n)$ with $w^k\in\Z^d$ for $k=1,\dots,n$,
and likewise for $l,u$, and where $b=(b^0,b^1,\dots,b^n)$ with $b^0\in\Z^r$
and $b^k\in\Z^s$ for $k=1,\dots,n$. Call an $n$-fold integer program {\em huge}
if $n$ is encoded in {\em binary}. More precisely, we are now given $t$
{\em types} of bricks, where each type $k=1,\dots,t$ has its cost
$w^k\in\Z^d$, lower and upper bounds $l^k,u^k\in\Z^d$, and
right-hand side $b^k\in\Z^s$. Also given are $b^0\in\Z^r$
and positive integers $n_1,\dots,n_t,n$ with $n_1+\cdots+n_t=n$, all encoded in binary.
A feasible point $x=(x^1,\dots,x^n)$ now must have first $n_1$
bricks of type $1$, next $n_2$ bricks of type $2$, and so on, with last
$n_t$ bricks of type $t$. Standard $n$-fold integer programming occurs as the
special case of $t=n$ and $n_1=\cdots=n_t=1$, and {\em symmetric} $n$-fold
integer programming occurs as the special case of $t=1$.

For $k=1,\dots,t$ the set of all possible bricks of type $k$ is the following,
$$S^k\ :=\ \{z\in\Z^d\ :\ A_2z=b^k\,,\ l^k\leq z\leq u^k\}\ .$$
We assume for simplicity that $S^k$ is finite for all $k$, which is the
case in most applications, such as in multiway table problems.
Let $\lambda^k:=(\lambda^k_z:z\in S^k)$ be a nonnegative integer tuple with entries
indexed by points of $S^k$. Each feasible point $x=(x^1,\dots,x^n)$ gives rise to
$\lambda^1,\dots,\lambda^t$ satisfying $\sum\{\lambda^k_z:z\in S^k\}=n_k$,
where $\lambda^k_z$ is the number of bricks of $x$ of type $k$ which are equal to $z$.
Let the {\em support} of $\lambda^k$ be
$\supp(\lambda^k):=\{z\in S^k:\lambda^k_z\neq 0\}$.
Then a {\em compact presentation of $x$} consists of the restrictions of $\lambda^k$
to $\supp(\lambda^k)$ for all $k$. However, the cardinality of $S^k$ may be
exponential in the binary encoding of the data $b^k,l^k,u^k$, so off hand
this presentation might be exponential as well. Nonetheless, we show the following.

\bt{N_Fold}
Consider data for the huge $n$-fold integer programming problem with
$t$ types over $(r,s)\times d$ bimatrix $A$, with $r,s,d$ fixed,
and with $w^k,b^k,l^k,u^k,n_1,\dots,n_t,n$ all encoded in binary.
Then the following three statements hold:
\begin{enumerate}
\item
If the problem is feasible then there is an optimal solution
which admits a compact presentation $\lambda^1,\dots,\lambda^t$
satisfying $|\supp(\lambda^k)|\leq 2^d$ for $k=1,\dots,t$.
\item
For $t$ fixed, the problem can be solved in polynomial time even
if the bimatrix $A$ is a variable part of the input and encoded in binary.
\item
For $A$ fixed and $t$ variable and encoded in unary, the {\em augmentation problem}
can be solved in polynomial time, namely, given a feasible point presented
compactly, we can either assert that it is optimal or find a better feasible point.
\end{enumerate}
\et

Here are some concluding remarks. First, it was shown in \cite{DO2} that {\em every}
bounded integer program can be isomorphically represented in polynomial time
for some $m$ and $n$ as some $3\times m\times n$ table problem.
So, by the above results, for any fixed $m$ we can handle integer programs with huge $n$.
Second, the results on threeway tables with line-sums can be extended to tables of
arbitrary fixed dimension and margins of any dimension. (A {\em $k$-margin} of a $d$-way
table is the sum of entries in some $(d-k)$-way subtable.)
We have the following theorem, stated without proof.

\bt{MultiwayTables}
Consider the huge multiway table problem over $m_1\times\cdots\times m_k\times n$
tables with $t$ types, with given margins of any dimension,
with $k,m_1,\dots,m_k$ fixed. It is in P
for every fixed $t$ and in NP intersect coNP for $t$ variable and unary encoded.
\et

\section{Proofs}

We begin by proving the three parts of Theorem \ref{N_Fold} one by one.
First, note that point $x=(x^1,\dots,x^n)$ is feasible in the huge
$n$-fold integer program only if each brick $x^i$ lies in some
$S^k=\{z\in\Z^d:A_2z=b^k,l^k\leq z\leq u^k\}$,
and $A_1\sum_{i=1}^n x^i=b^0$. So our assumption that each $S^k$
is finite implies that the set of feasible points is finite as well.
Therefore, if the program is feasible then it has an optimal solution.

The proof of part (1) makes use of a nice argument of Eisenbrand-Shmonin \cite{ES}.

\vskip.2cm{\em Proof of Theorem \ref{N_Fold} part (1).}
Suppose the huge $n$-fold program is feasible. Then, as explained above,
there is an optimal solution. Let $x$ be an optimal solution with minimum value
$\sum_{i=1}^n\|x^i\|^2$ with $\|z\|^2=\sum_{j=1}^dz_j^2$ the Euclidean norm squared.
Let $\lambda^1,\dots,\lambda^t$ be a compact presentation of $x$.
Suppose indirectly that we have $|\supp(\lambda^j)|>2^d$. Then there are
two vectors $y'\neq y''$ in $\supp(\lambda^j)$ having the same parity on each coordinate,
implying $y:={1\over2}(y'+y'')\in S^j$. For $k=1,\dots, t$ define $\mu^k$ on
$S^k$ to be the same as $\lambda^k$ except that $\mu^j_{y'}:=\lambda^j_{y'}-1$,
$\mu^j_{y''}:=\lambda^j_{y''}-1$ and $\mu^j_y:=\lambda^j_y+2$. Let ${\bar x}$ be
the vector whose compact presentation is given by the $\mu^k$. Then
$$\sum_{i=1}^n x^i-\sum_{i=1}^n {\bar x}^i\ =\
\sum_{k=1}^t\sum_{z\in S^k}\lambda^k_z z-\sum_{k=1}^t\sum_{z\in S^k}\mu^k_z z
\ =\ y'+y''-2y\ =\ 0$$
and therefore $A_1\sum_{i=1}^n {\bar x}^i=A_1\sum_{i=1}^n x^i=b^0$ so ${\bar x}$ is
also feasible. Furthermore,
$$wx-w{\bar x}\ =\ \sum_{k=1}^t w^k\sum_{z\in S^k}(\lambda^k_z-\mu^k_z) z\ =\
w^j(y'+y''-2y)\ =\ 0$$
and therefore ${\bar x}$ is also optimal. But now we have
\begin{eqnarray*}
\sum_{i=1}^n \|x^i\|^2-\sum_{i=1}^n \|{\bar x}^i\|^2 & = &
\sum_{k=1}^t\sum_{z\in S^k}(\lambda^k_z-\mu^k_z)\|z\|^2
\ =\ \|y'\|^2+\|y''\|^2-2\|y\|^2 \\
&=&\|y'\|^2+\|y''\|^2-2\|{1\over2}(y'+y'')\|^2\ =\ {1\over2}\|y'-y''\|^2\ >\ 0
\end{eqnarray*}
which is a contradiction to the choice of $x$. This completes the proof.
\epr

\vskip.5cm
The proof of part (2) uses the following beautiful
result of \cite{GR} building on \cite{ES}.

\bp{IntegereCones}{\bf (Goemans-Rothvo\ss).}
Fix $d,t$. Let $S^k=\{z\in\Z^d:A^k z\leq a^k\}$ be finite for $k=1,\dots,t$,
and let $T=\{z\in\Z^d:Bz\leq b\}$. Then, in polynomial time, we can
decide if there are nonnegative integer tuples $\lambda^k=(\lambda^k_z:z\in S^k)$
such that $\sum_{k=1}^t\sum_{z\in S^k}\lambda^k_z z\in T$,
in which case we can compute such $\lambda^k$ of polynomial support.
\ep

\vskip.2cm{\em Proof of Theorem \ref{N_Fold} part (2).}
We make use of points in $z\in\Z^{t(d+1)}$ and index each such point by
$z=(z^1_0,z^1,\dots,z^t_0,z^t)$ with $z^k_0\in\Z$ and $z^k\in\Z^d$ for $k=1,\dots,t$.
Let $L\leq U$ be two integers. Define the following sets
$S^1,\dots,S^t$ and $T$ in $\Z^{t(d+1)}$,
$$S^k\ :=\ \left\{
z\in\Z^{t(d+1)}\ :\ z^k_0=1,\ A_2 z^k=b^k,\ l^k\leq z^k\leq u^k,\
z^i_0=0,\ z^i=0,\ i\neq k
\right\}\ ,$$
$$T\ :=\ \left\{
y\in\Z^{t(d+1)}\ :\ y^1_0=n_1,\dots,y^t_0=n_t,\ A_1\sum_{k=1}^t y^k=b^0,\
L\leq \sum_{k=1}^t w^k y^k\leq U
\right\}\ .$$
Now suppose that $x=(x^1,\dots,x^n)$ is a feasible point in the huge $n$-fold integer
program, with objective function value $wx=v$ which satisfies $L\leq v\leq U$.
Note that $\{z^k\in\Z^d:z\in S^k\}$ is the set of
possible bricks of $x$ of type $k$, and let $\lambda^k=(\lambda^k_z:z\in S^k)$
for $k=1,\dots,t$ be nonnegative integer tuples with $\lambda^k_z$
the number of bricks of $x$ of type $k$ which are equal to $z^k$.
Let $y:=\sum_{i=1}^t\sum_{z\in S^i}\lambda^i_z z\in \Z^{t(d+1)}$.

Since $x$ is feasible, we have
$$y^k_0\ =\ \sum_{i=1}^t\sum_{z\in S^i}\lambda^i_z z^k_0\ =\
\sum_{z\in S^k}\lambda^k_z \ =\ n_k\ ,\quad k=1,\dots,t\ ,$$
$$A_1\sum_{k=1}^t y^k\ =\ A_1\sum_{k=1}^t\sum_{i=1}^t\sum_{z\in S^i}\lambda^i_z z^k
\ =\ A_1\sum_{k=1}^t\sum_{z\in S^k}\lambda^k_z z^k\ =\ A_1\sum_{j=1}^n x^j\ =\ b^0\ ,$$
and
$$\sum_{k=1}^t w^ky^k\ =\ \sum_{k=1}^t w^k\sum_{z\in S^k}\lambda^k_z z^k\ =\
\sum_{k=1}^t w^k\sum\{x^j:x^j\ \mbox{has type}\ k\}\ =\ wx\ =\ v\ .$$
So $y$ is a nonnegative integer combination of points
of $\bigcup_{k=1}^t S^k$ which lies in $T$.

Conversely, suppose $y=\sum_{i=1}^t\sum_{z\in S^i}\lambda^i_z z$ is a nonnegative
integer combination of points of $\bigcup_{k=1}^t S^k$ and $y\in T$,
and let $v:=\sum_{k=1}^t w^ky^k$. Then
$$\sum_{z\in S^k}\lambda^k_z\ =\ \sum_{i=1}^t\sum_{z\in S^i}\lambda^i_z z^k_0\ =\
 y^k_0\ =\ n_k\ ,\quad k=1,\dots,t\ ,$$
so we can construct a vector $x=(x^1,\dots,x^n)$ with $\lambda^k_z$ bricks of type
$k$ which are equal to $z^k$ for $k=1,\dots,t$ and all $z\in S^k$. We then have
$$A_1\sum_{j=1}^n x^j\ =\ A_1\sum_{k=1}^t\sum_{z\in S^k}\lambda^k_z z^k\ =\
A_1\sum_{k=1}^t\sum_{i=1}^t\sum_{z\in S^i}\lambda^i_z z^k
\ =\ A_1\sum_{k=1}^t y^k\ =\ b^0\ ,$$
so $x$ is feasible in the huge $n$-fold program, and has objective function value
$$wx\ =\ \sum_{k=1}^t w^k\sum\{x^j:x^j\ \mbox{has type}\ k\}\ =\
\sum_{k=1}^t w^k\sum_{z\in S^k}\lambda^k_z z^k\ =\ \sum_{k=1}^t w^ky^k\ =\ v\ .$$

Since $d$ and $t$ are fixed and $S^1,\dots,S^t$ are finite, applying
Proposition \ref{IntegereCones} to the $S^k$ and $T$ in $\Z^{t(d+1)}$, we can in
polynomial time decide if there is a feasible point $x$ in the $n$-fold program
with objective function value in the interval $[L,U]$, and if there is, find a compact
presentation $\lambda^1,\dots,\lambda^t$ of $x$, each with polynomial support.

Now, using the algorithm for integer programming in fixed dimension \cite{Len},
we find $L_k:=\min\{w^kz:z\in S^k\}$ and $U_k:=\max\{w^kz:z\in S^k\}$ for $k=1,\dots,t$.
Then any feasible point in the $n$-fold program has objective value in the
interval $[\sum_{k=1}^t n_kL_k,\sum_{k=1}^t n_kU_k]$, and so
by binary search on that interval and repeated application of the above procedure
starting with $L:=\sum_{k=1}^t n_kL_k$ and $U:=\sum_{k=1}^t n_kU_k$,
we can solve the huge symmetric $n$-fold integer program in polynomial time.
\epr

For the proof of part (3) we need to review some facts about Graver bases.
We introduce a partial order $\sqsubseteq$ on $\Z^n$ by $x\sqsubseteq y$ if
$x_iy_i\geq 0$ and $|x_i|\leq |y_i|$ for $i=1,\ldots,n$. The {\em Graver basis}
of an integer $m\times n$ matrix $B$ is the set $\G(B)\subset\Z^n$
of all $\sqsubseteq$-minimal elements in $\{x\in\Z^n:Bx=0,x\neq 0\}$.
It is well known that the Graver basis is a {\em test set} for any integer program
of the form $\min\{wx:x\in\Z^n,Bx=b,l\leq x\leq u\}$ defined by $B$, that is,
if $x$ is a feasible but not optimal in that program, then there is an
element $y\in\G(B)$ such that $x+y$ is feasible and better, see \cite{Onn}.
It is also known that the Graver basis of any integer matrix
is finite, but it may be exponentially large. However, Graver bases of
$n$-fold products are well behaved as we now explain. Let $n\geq g$ and let
$\G(A^{(g)})$ be the Graver basis of the $g$-fold product of a bimatrix $A$.
An {\em $n$-lifting} of $h=(h^1,\dots,h^g)\in\G(A^{(g)})$ is any vector
$y=(y^1,\dots,y^n)$ such that for some $1\leq i_1<\cdots<i_g\leq n$,
$y^{i_1}=h^1,\dots,y^{i_{g}}=h^g$ and all other bricks of $y$ are $0$.

\bp{Graver_bases}{\bf (see \cite{Onn}).}
For every bimatrix $A$ there is a constant $g(A)\in\Z_+$, called the
{\em Graver complexity} of $A$, such that for all $n\geq g(A)$, the Graver basis $\G(A^{(n)})$
consists precisely of all $n$-liftings of elements of the Graver basis $\G(A^{(g(A))})$.
\ep

\vskip.2cm{\em Proof of Theorem \ref{N_Fold} part (3).}
Let $g:=g(A)$ be the Graver complexity of $A$. If $n<g$ then we can solve
the $n$-fold program, and in particular the augmentation problem, using integer
programming in fixed dimension $nd<gd$ in polynomial time \cite{Len}.
So assume $n\geq g$. Let $\lambda^1,\dots,\lambda^t$ be a compact presentation
of a feasible point $x$. Suppose $x$ is not optimal. We show how to find
$y\in\G(A^{(n)})$ such that $x+y$ is feasible and better.
In fact, we can find $y\in\G(A^{(n)})$ and step size $\alpha\in\Z_+$ such that
$x+\alpha y$ is feasible and attains the best possible improvement attainable
by any multiple of any Graver basis element.
Consider any $h=(h^1,\dots,h^g)\in\G(A^{(g)})$.
Let $\biguplus_{k=1}^t \supp(\lambda^k)$ be the disjoint union of the supports
of the $\lambda^k$ (so a point which happens to be in the support of more
than one $\lambda^k$ appears more than once). Consider a mapping
$$\phi:\{h^1,\dots,h^g\}\rightarrow \biguplus_{k=1}^t \supp(\lambda^k)
\ : \ h^i\mapsto z^i\in\supp(\lambda^{k(i)})\ .$$
Such a mapping provides a compact way of prescribing an $n$-lifting $y$ of $h$.
For such a lifting and any $\alpha\in\Z_+$, we will have that $x+\alpha y$ is feasible
and better than $x$ if the following conditions hold:
(1) $|\phi^{-1}(z)|\leq \lambda^k_z$ for $k=1,\dots,t$ and all $z\in\supp(\lambda^k)$;
(2) $l^{k(i)}\leq z^i+\alpha h^i\leq u^{k(i)}$ for $i=1,\dots,g$;
(3) $\sum_{i=1}^g w^{k(i)}\alpha h^i<0$.
(Note that each $h^i$ satisfies $A_2h^i=0$ and hence
$A_2(z^i+\alpha h^i)=A_2z^i=b^{k(i)}$ holds automatically.)
Now, it can be checked if these conditions hold, say, with $\alpha=1$, and if they do,
the maximum $\alpha$ for which they hold be computed, easily in polynomial time.
Moreover, a compact presentation $\mu^1,\dots,\mu^t$ of the new better point
$x+\alpha y$ can be obtained as follows. Begin by defining $\mu^k:=\lambda^k$
for $k=1,\dots,t$. Now, for $i=1,\dots,g$, set
$$\mu^{k(i)}_{z^i}\ :=\ \mu^{k(i)}_{z^i}-1,\quad
\mu^{k(i)}_{z^i+\alpha h^i}:=\mu^{k(i)}_{z^i+\alpha h^i}+1\ .$$
This provides a compact presentation of the new
feasible and better point $x+\alpha y$.

Now, since the bimatrix $A$ is fixed, so is its Graver complexity $g=g(A)$ and
hence so is the number of elements $h\in\G(A^{(g)})$. Moreover, the number
of possible lifting mappings $\phi$ of $h$ is $|\biguplus_{k=1}^t \supp(\lambda^k)|^g$
which is polynomial in the size of the input which includes the compact presentation
$\lambda^1,\dots,\lambda^t$ of $x$. So by going over all $h\in\G(A^{(g)})$ and
$\phi$ we can either find that there is no feasible better point of the form
$x+y$ and conclude that $x$ is optimal, or find $h\in\G(A^{(g)})$,
mapping $\phi$, $\alpha\in\Z_+$, and compact presentation $\mu^1,\dots,\mu^t$
of that $x+\alpha y$ which gives best improvement.
\epr

We proceed to establish Theorem \ref{ThreeWayLineSum}.

\vskip.2cm{\em Proof of Theorem \ref{ThreeWayLineSum}.}
The huge threeway table problem can be formulated as a huge $n$-fold integer
programming problem as follows. Let $r=d=lm$ and $s=l+m$, and let $A_1=I_{lm}$
be the $lm\times lm$ identity matrix and $A_2$ be the $(l+m)\times lm$ incidence
matrix of the complete bipartite graph $K_{l,m}$. (So $A_2$ is itself an
$m$-fold product $A=B^{(m)}$ with $B$ the bimatrix having $B_1=I_l$ and $B_2$
a row of $l$ ones.) Index $l\times m\times n$ tables as $x=(x^1,\dots,x^n)$
with $x^k=(x_{1,1},\dots,x_{l,1},\dots,x_{1,m},\dots,x_{l,m})$,
set all costs $w^k_{i,j}:=0$, all lower bounds $l^k_{i,j}:=0$,
and all upper bounds $u^k_{i,j}:=\infty$, and arrange the row-sums, column-sums,
and line sums suitably in the right-hand side vector $b$, with $b^0=g$ and
$b^k=(f^k,e^k)$ for $k=1,\dots,t$. This encodes the huge table problem
as a huge $n$-fold integer program with a {\em fixed} bimatrix $A$.

Now, the statement of Theorem \ref{ThreeWayLineSum} for fixed $t$ follows from
Theorem \ref{N_Fold} part (2), since solving the optimization problem
in particular enables to decide feasibility,

Next, we prove the statement of Theorem \ref{ThreeWayLineSum} for variable $t$.
First note that if the problem is feasible then, by Theorem \ref{N_Fold}
part (1) it has an optimal solution with compact presentation $\lambda^1,\dots,\lambda^t$
satisfying $|\supp(\lambda^k)|\leq 2^{lm}$ for $k=1,\dots,t$, providing a polynomial
certificate for feasibility and showing the problem is in NP.

We proceed to show that the problem is in coNP.
Let $c:=2r+s=2lm+(l+m)$ and define an $(r,s)\times c$ bimatrix $C$ with blocks
$C_1:=[A_1,\ I_{lm},\ 0_{lm\times (l+m)}]$ and $C_2:=[A_2,\ 0_{(l+m)\times lm},\ I_{l+m}]$
with $A_1,A_2$ as above. Define a huge $n$-fold program over the bimatrix $C$
as follows. As before, set all lower bounds to $0$ and all upper bounds to $\infty$.
Without loss of generality assume that all row-sums, column-sums, and line-sums
are nonnegative, else there is no feasible table, and arrange them in the right-hand
side vector $b$ as before. Now index the variables as $v=(v^1,\dots,v^n)$ with each
brick of the form $v^i=[x^i,y^i,z^i]$ with $x^i,y^i\in\Z_+^{lm}$ and $z^i\in\Z_+^{l+m}$,
with $x^i$ interpreted as a layer of the sought table and $y^i,z^i$ as slacks.
There are again $t$ types where a brick $v$ of each type $k=1,\dots,t$ must satisfy
$C_2v=b^k$, and the multiplicities are $n_1,\dots,n_t$ with
$n_1+\cdots+n_t=n$ as given for the table problem.

Now, for this auxiliary program we can always write down a compact presentation of a
feasible point $v$ defined as follows. We use the brick $a^0:=[0,b^0,b^1]$ of type $1$
with multiplicity $\lambda^1_{a^0}:=1$, the brick $a^1:=[0,0,b^1]$ of type $1$ with
multiplicity $\lambda^1_{a^1}:=n_1-1$, and for $k=2,\dots,t$, the brick $a^k:=[0,0,b^k]$
of type $k$ with multiplicity $\lambda^k_{a^k}:=n_k$. Clearly, these bricks are
nonnegative and each brick of type $k$ indeed satisfies $C_2v=b^k$. Moreover,
$C_1\sum_{i=1}^n v^i=C_1\left(a^0+\sum_{k=1}^t\lambda^k_{a^k}a^k\right)=b^0$.

Now, we consider the problem of minimizing the sum of all slack variables.
Note that the value of this sum will be always nonnegative, and will be $0$ if and only
if all slacks are $0$, which holds if and only if the restriction $x:=(x^1,\dots,x^n)$
of $v$ is a feasible table, which holds if and only if the table problem is feasible.

Suppose now the table problem is infeasible. By Theorem \ref{N_Fold}
part (1), the auxiliary program has an optimal solution, i.e. $v$ minimizing the sum
of slack variables, with compact presentation $\lambda^k$ satisfying
$|\supp(\lambda^k)|\leq 2^{2lm+l+m}$ for $k=1,\dots,t$. Now, using this compact presentation,
we can compute the sum of slacks and verify that it is positive, and using
Theorem \ref{N_Fold} part (3), we can verify that $v$ is indeed an optimal
solution, in polynomial time. This proves that the problem is in coNP.
\epr


\begin{thebibliography}{}

\bibitem{AT}
Aoki, S., Takemura, A.:
Minimal basis for connected Markov chain over $3\times3\times K$
contingency tables with fixed two-dimensional marginals.
Australian and New Zealand Journal of Statistics 45:229--249 (2003)

\bibitem{DHOW}
De Loera, J., Hemmecke, R., Onn, S., Weismantel, R.:
N-fold integer programming.
Discrete Optimization 5:231--241 (2008)

\bibitem{DO1}
De Loera, J., Onn, S.:
The complexity of three-way statistical tables.
SIAM Journal on Computing 33:819--836 (2004)

\bibitem{DO2}
De Loera, J., Onn, S.:
All linear and integer programs are slim 3-way transportation programs.
SIAM Journal on Optimization 17:806--821 (2006)

\bibitem{ES}
Eisenbrand, F., Shmonin, G.:
Carath\'eodory bounds for integer cones.
Operations Research Letters 34:564--568 (2006)

\bibitem{GR}
Goemans, M.X., Rothvo\ss, T.:
Polynomiality for Bin Packing with a Constant Number of Item Types.
Symposium on Discrete Algorithms 25:830--839 (2014)

\bibitem{HOR}
Hemmecke, R., Onn, S., Romanchuk, L.:
N-fold integer programming in cubic time.
Mathematical Programming 137:325--341 (2013)

\bibitem{HOW}
Hemmecke, R., Onn, S., Weismantel, R.:
A polynomial oracle-time algorithm for convex integer minimization.
Mathematical Programming 126:97--117 (2011)

\bibitem{HS}
Ho\c sten, S., Sullivant, S.:
Finiteness theorems for Markov bases of hierarchical models.
Journal of Combinatorial Theory Series A 114:311--321 (2007)

\bibitem{Len}
Lenstra, H.W., Jr.:
Integer programming with a fixed number of variables.
Mathematics of Operations Research 8:538--548 (1983)

\bibitem{Onn}
Onn, S.: Nonlinear Discrete Optimization.
Zurich Lectures in Advanced Mathematics,
European Mathematical Society (2010),
available online at: {\tt http://ie.technion.ac.il/$\sim$onn/Book/NDO.pdf}

\bibitem{SS}
Santos, F., Sturmfels, B.:
Higher Lawrence configurations.
Journal of Combinatorial Theory Series A 103:151--164 (2003)

\end{thebibliography}
\end{document}